\documentclass[reqno]{amsart}
\usepackage{amssymb,amsthm,amscd}
\usepackage[pdftex]{graphics}

\newcommand{\fp}{{\mathfrak p}}
\newcommand{\fq}{{\mathfrak q}}

\newcommand{\Q}{{\mathbb Q}}

\newcommand{\fa}{\mathfrak a}
\newcommand{\fb}{\mathfrak b}
\newcommand{\fc}{\mathfrak c}
\newcommand{\fd}{\mathfrak d}

\newcommand{\fthr}{\mathfrak 3}

\newcommand{\cO}{{\mathcal O}}

\newcommand{\Cl}{{\operatorname{Cl}}}

\newcommand{\lra}{\longrightarrow}

\newcommand{\Kts}{K^{\times\,2}}

\newfont{\cyr}{wncyb10}
\newcommand{\TS}{\mbox{\cyr Sh}}

\newcommand{\num}{\operatorname{num}}
\newcommand{\den}{\operatorname{den}}

\newtheorem{thm}{Theorem}
\newtheorem{prop}[thm]{Proposition}
\newtheorem{lem}[thm]{Lemma}
\newtheorem{cor}[thm]{Corollary}

\title{Binomial Squares in Pure Cubic Number Fields}

\author{F. Lemmermeyer}
\email{hb3@ix.urz.uni-heidelberg.de}
\address{M\"orikeweg 1, 73489 Jagstzell, Germany}

\begin{document}

\begin{abstract}
Let $K = \Q(\omega)$ with $\omega^3 = m$ be a pure cubic 
number field. We show that the elements $\alpha \in K^\times$ whose 
squares have the form $a - \omega$ form a group isomorphic to the 
group of rational points on the elliptic curve $E_m: y^2 = x^3 - m$.
This result will allow us to construct unramified quadratic extensions
of pure cubic number fields $K$ without any calculations in $K$.
\end{abstract}

\maketitle

\section*{Introduction}

The first attempt of solving the diophantine equation $y^2 = x^3 + 1$ 
(or, equivalently, determining the integral points on this elliptic 
curve) probably is writing it in the form $(y-1)(y+1) = x^3$ and 
using unique factorization. If $y$ is odd, one possibility we have
to consider is $y-1 = 2a^3$ and $y+1 = 4b^3$ for integers $a, b$. 
This implies $a^3 - 2b^3 = -1$, which has the obvious solution
$(a,b) = (1,1)$, giving the solution $(x,y) = (2,3)$ of the 
original equation. Showing that $a^3 - 2b^3 = -1$ does not have
any other solution is a slightly technical task. A famous result
due to Delaunay and Nagell shows that equations $a^3 - mb^3 = 1$
for noncubes $m$ have at most one integral solution. For proving
this result one needs to study units of the form $a - b \sqrt[3]{m}$
in pure cubic number fields and show that these units, with a few
exceptions, cannot be powers of other units.

In this article we shall investigate squares of the form 
$a + b \sqrt[3]{m}$ in pure cubic number fields and explain
why their occurrence in diophantine problems related to elliptic
curves is quite natural. We will also show how to apply our results
to the construction of unramified quadratic extensions of pure cubic
number fields; in particular, we will give these extensions for all
values $2 \le m \le 113$ for which $K = \Q(\sqrt[3]{m}\,)$ has even
class number.

Let us start by giving a few nontrivial\footnote{Trivial examples are 
$r^2 = r^2 + 0 \sqrt[3]{m}$ and ${\sqrt[3]{m^2}\,}^2 = 0 + m \sqrt[3]{m}$.} 
binomial squares of the form $a - \sqrt[3]{2}$ in $\Q(\sqrt[3]{2}\,)$:
\begin{align*}
 (1 - \sqrt[3]{2} - \sqrt[3]{4}\,)^2 & = 5 - \sqrt[3]{4}, \\
(9 - 6 \sqrt[3]{2} - 2 \sqrt[3]{4}\,)^2 & = 129 - 100  \sqrt[3]{2}, \\
 (16641 - 25800 \sqrt[3]{2} - 20000 \sqrt[3]{4}\,)^2
    & = 2340922881 - 58675600 \sqrt[3]{2}, 
\end{align*}
where $58675600 = 7660^2$.

This abundance of examples in $\Q(\sqrt[3]{2}\,)$ should not mislead
the readers into thinking that this is a typical phenomenon; in fact,
there are no nontrivial squares of the form $a - \sqrt[3]{3}$ at all 
in $\Q(\sqrt[3]{3}\,)$.

\section{The Group Law}\label{S1}

Fix a cubefree integer $m$, let $K = \Q(\sqrt[3]{m}\,)$ denote the
corresponding pure cubic number field, and consider the set 
$$ S_m = \{\alpha \in K^\times : \alpha^2 = a - \sqrt[3]{m}: a \in \Q^\times\}. $$
Writing $\alpha = r + s \omega + t \omega^2$ with $\omega = \sqrt[3]{m}$
we find that the condition $\alpha^2 = a - \omega$ is equivalent to the 
system of equations
\begin{align}
\label{E1}  2rt + s^2  & = 0, \\
\label{E2}  2rs + mt^2 & = -1, \\
\label{E3}  2mst + r^2 & = a. 
\end{align}

Since $t = 0$ implies $s = 0$ and $a = r^2$ (which is the trivial solution), 
we may assume that $t \ne 0$; solving (\ref{E1}) for $r$ and plugging the 
resulting equation
$r = - s^2/2t$ into (\ref{E2}) we find $-s^3/t + mt^2 = -1$ which,
after dividing through by $-t^2$ gives the point
$$ P_\alpha = (x,y) = \Big(\frac st, \frac 1t \Big)  
     \quad \text{ on the elliptic curve } y^2 = x^3 - m. $$
The missing parameters $a$ and $r$ are given by $r/t = -\frac12 (s/t)^2$
and by (\ref{E3}).

Conversely, every affine point $(x,y) \in E_m(\Q)$ gives a unique element 
of $S_m$ via
\begin{equation}\label{Ebij}
 t = \frac1y, \quad s = \frac xy, \quad r = \frac{x^2}{2y}, \quad 
   \text{and} \quad a = \frac{x^4 + 8mx}{4y^2}. 
\end{equation}
Of course we make the point at infinity on $E_m(\Q)$  correspond to the
(class of the) trivial element $1 \in S_m$.

We have proved

\begin{thm}\label{Tm}
Let $m$ be an integer that is not a cube, let $K = \Q(\omega)$ with 
$\omega^3 = m$ denote the corresponding pure cubic field, and let 
$E_m$ denote the elliptic curve $y^2 = x^3 - m$. There is a bijection 
between the rational points $(x,y) \in E_m(\Q)$ and the elements 
$\alpha  \in K^\times$ with $\alpha^2 = a - \omega$ for $a \in \Q^\times$. 
In fact, if $\alpha = r + s \omega + t \omega^2$ satisfies 
$\alpha^2 = a - \omega$, then $(x,y) \in E_m(\Q)$ for $x = \frac st$ 
and $y = \frac1t$. Clearly $\alpha = 1$ corresponds to the point at 
infinity on $E$, and multiplication by $-1$ in $K^\times$ corresponds to
multiplication by $-1$ in $E_m(\Q)$.

Conversely, an affine point $(x,y)$ corresponds to 
$$ \alpha = -\frac{x^2}{2y} + \frac xy \omega + \frac1y \omega^2 \quad
    \text{ with } \alpha^2 = 
        \frac{x^4 + 8mx}{4y^2}  -  \omega. $$
\end{thm}

Multiplying through by $4y^2$ gives the simpler identity
$$ (x^2 - 2x \omega - 2 \omega^2)^2 = x^4 + 8mx - 4y^2 \omega, $$
which can easily be verified directly.

\medskip\noindent{\bf Remark.}
If $\alpha = r + s \omega + t \omega^2$ with $\omega = \sqrt[3]{m}$
satisfies $\alpha^2 = a - \omega$ for some $a \in \Q^\times$, then 
$\alpha_1 = r - s \omega + t \omega^2$ has the property that 
$\alpha_1^2 = c - d \omega$ for rational numbers $c$, $d$; this
is due to the fact that Eqn. (\ref{E1}) remains invariant under 
$s \to -s$.

\begin{cor}
If $E_m: y^2 = x^3 - m$ has no rational point except the point at infinity, 
then the only squares of the form $a - \omega$ in $K = \Q(\omega)$ with 
$a \in \Q^\times$ and $\omega = \sqrt[3]{m}$ are those with $a = r^2$ 
for $r \in \Q^\times$.
\end{cor}

\medskip\noindent{\bf Remark.} Torsion points of order $3$ on $E_m$ 
may be neglected; in fact, the torsion points $(0, \pm k)$ on 
$y^2 = x^3 + k^2$ give rise to the trivial solutions
$ \sqrt[3]{k^2}^2 = k \sqrt[3]{k}. $
\medskip

Some expressions occurring in Thm. \ref{Tm} have a natural
explanation in terms of the group law on elliptic curves.

 The group law on $E_m$ is given 
by the following formulas: given rational points $(x_1,y_1)$ and 
$(x_2,y_2)$, set 
$$ \lambda = \begin{cases} \smallskip
                \frac{y_2-y_1}{x_2-x_1}   & \text{ if } x_2 \ne x_1, \\
                \ \frac{3x_1^2}{2y_1}     & \text{ if } x_2  =  x_1.
              \end{cases} $$
Then $(x_1,y_1) + (x_2,y_2) = (x_3,y_3)$ with
$$ x_3 = \lambda^2 - x_1 - x_2, \quad
   y_3 = \lambda x_3 + y_1 - \lambda x_1. $$
The duplication formula for $P = (x,y)$ is given by
$$ 2P = \Big(\frac{x^4 + 8mx}{4y^2}, 
                   - \frac{x^6 - 20mx^3 - 8m^2}{8y^3} \Big). $$
Comparing this formula with (\ref{Ebij}) we immediately get 

\begin{cor}\label{C2}
The element $a - \omega$ is a square in $K$ if and only if
$a = x_{2P}$ for some point $P \in E_m(\Q)$. In this case,
$N(a - \omega) = y_{2P}^2$.
\end{cor}

The last claim is a simple calculation:
$$ \Big(\frac{x^4 + 8mx}{4(x^3-m)} \Big)^3 - m = 
                \Big(\frac{x^6 - 20mx^3 - 8m^2}{8y^3}\Big)^2. $$

This result can be used to simplify a technical point in \cite{Lem}.
There we have proved (a special case of) the following result via
descent:

\begin{prop}
Let $m$ be a cubefree integer, and let $E: y^2 = x^3 - m$ denote an 
elliptic curve. If $P \in E(\Q)\setminus 2E(\Q)$ is written in the
form $P = (\frac{r}{t^2},\frac{s}{t^3})$ for integers $r, s, t$ with 
$\gcd(r,t) = \gcd(s,t) = 1$, then $\alpha = r - t^2 \omega \in K^\times$
is not a square in $K$.
\end{prop}

Cor. \ref{C2} shows immediately that $\alpha$ is a square in $K$ if 
and only if $x_P = \frac{r}{t^2} = x_{2Q}$ is the $x$-coordinate of a 
point $2Q$ on $E_m$.

\medskip \noindent {\bf Example.}
Consider the element $\beta = 5 - \sqrt[3]{4}$ in $K = \Q(\sqrt[3]{2}\,)$, 
whose norm $5^3 - 4 = 11^2$ is a square. For deciding whether $\beta$ is a 
square in $K$ we observe that the point $P = (5,11)$ on 
$y^2 = x^3 - 4$ is a multiple of $2$ since $P = 2 (2,-2)$. Thus 
$\beta = \alpha^2$ is a square, and the formulas above give 
$\alpha = -1 + \sqrt[3]{4} + \frac12 \sqrt[3]{16}
       = -1 + \sqrt[3]{2} + \sqrt[3]{4}$ 
as well as $N\beta = N(\alpha)^2 = 11^2$.

\medskip\noindent {\bf The Example of Bachet-Fermat.}
Consider the curve $E_2: y^2 = x^3 - 2$. The integral point 
$(x,y) = (3,5)$ corresponds to 
$$ \alpha_P = - \frac9{10} + \frac 35 \sqrt[3]{2} + \frac15 \sqrt[3]{4} \quad
   \text{with} \quad
   \alpha_P^2 = \frac{129}{100} - \sqrt[3]{2}. $$
Multiplying through by $10^2$ shows that $P$ gives rise to
$$ (-9 + 6 \sqrt[3]{2} + 2 \sqrt[3]{4}\,)^2 = 129 - 100  \sqrt[3]{2}. $$

Observe that $2P = (\frac{129}{100}, \frac{383}{1000})$ corresponds to 
$$ \alpha_{2P} = -\frac{16641}{7660} + \frac{1290}{383} \omega 
                   + \frac{1000}{383} \omega^2 $$
with $\alpha_{2P}^2 = \frac{2340922881}{58675600} - \omega$. Finally
we remark that 
\begin{align*}
  3P & = \Big(\frac{164323}{171^2}, - \frac{66234835}{171^3} \Big), \\
  \alpha_{3P} & = 
           \frac{27002048329}{22652313570} 
         - \frac{28099233}{66234835} \sqrt[3]{2}
         - \frac{5000211}{66234835} \sqrt[3]{4}. \end{align*}

\section{Explicit Multiplication Formulas}

Recall that we have constructed a bijection between elements 
$\alpha \in K^\times$ with $\alpha^2 = a - \omega$ and rational
points on the elliptic curve $E_m: y^2 = x^3 - m$. The group 
structure on $E_m(\Q)$ induces, by transport of structure, a 
group law on $S_m$. In this section we will give explicit
formulas for the group law on $S_m$.

To this end assume that $\alpha_j = r_j + s_j \omega + t_j \omega^2$ 
($j = 1, 2$) are elements whose squares have the form $a_j - \omega$. 
Then these elements correspond to the points
$$ P_j = \Big(\frac{s_j}{t_j}, \frac1{t_j}\Big)  $$
on the elliptic curve $E_m: y^2 = x^3 - m$. The sum $P_3 = P_1 + P_2$
corresponds to an element $\alpha_3 \in K^\times$ with 
$\alpha_3^2 = a_3 - \omega$, and we can compute formulas for 
$\alpha_3 = \alpha_1 * \alpha_2$ by using the group law on 
elliptic curves. The ``multiplication formulas'' for the
$\alpha_j$ have little if anything to do with multiplication
in $K^\times$ and are rather complicated. It remains to be seen
whether there is any geometric interpretation of these formulas.

If $\frac{s_1}{t_1} \ne \frac{s_2}{t_2}$, then we find
\begin{align*}
 x_3 & = \bigg( \frac{\frac1{t_2} - \frac1{t_1}}
                    {\frac{s_2}{t_2} - \frac{s_1}{t_1}} \bigg)^2 
           - \frac{s_1}{t_1} - \frac{s_2}{t_2}
       = \Big( \frac{t_1 - t_2}{t_1s_2 - t_2s_1} \Big)^2  
                  - \frac{s_1}{t_1} - \frac{s_2}{t_2}, \\
 y_3 & =  \frac{t_1 - t_2}{t_1s_2 - t_2s_1}\Big(x_3 - \frac{s_1}{t_1}\Big) 
            + \frac1{t_1}. 
\end{align*}
From these values we can compute the coefficients of
$$ \alpha_1 * \alpha_2 = \alpha_3 = r_3 + s_3 \omega + t_3 \omega^2. $$
For what it's worth, the corresponding values of $s_3$ and $t_3$ are given by
$$ s_3 = \frac{x_3}{y_3} = \frac{\num(s)}{\den(s)}, \qquad
   t_3 = \frac1{y_3} =  \frac{\num(t)}{\den(t)} $$
with 
\begin{align*}
 \num(s) & = (s_1t_2 - s_2t_1)((s_1t_2+s_2t_1)(s_1t_2-s_2t_1)^2
                - t_1t_2(t_1-t_2)^2) \\
 \den(s) & = t_1t_2(t_1-t_2)^3 + (s_1t_2-s_2t_1)^2
             (s_1t_2^2 + 2(s_2 - s_1)t_2t_1 - s_2t_1^2) \\
 \num(t) & = (s_2t_1 - s_1t_2)^3 t_1t_2, \\
 \den(t) & = t_1t_2 (t_1-t_2)^3 + (s_1t_2 - s_2t_1)^2
                   (s_1t_2^2 + 2(s_2-s_1)t_1t_2 - s_2t_1^2). 
\end{align*}
Setting 
$$ S_- = s_1t_2 - s_2t_1, \quad S_+ = s_1t_2+s_2t_1,
   T_- = t_1-t_2 \quad \text{and} \quad T_+ = t_1t_2, $$
as well as
$$ \Sigma = s_1t_2^2 + 2(s_2-s_1)t_1t_2 - s_2t_1^2
          = (s_2-s_1)T_+ - S_+T_-, $$
the multiplication formulas become
$$ s_3 = \frac{S_-^3S_+ - S_-T_-^2T_+}{T_-^3T_+ + S_-^2 \Sigma}, \qquad
  t_3 = - \frac{S_-^3T_+}{T_-^3T_+ + S_-^2\Sigma}. $$ 
Let us check these formulas by ``multiplying'' the elements
$$ \alpha_1 =   \frac{9}{10} - \frac{3}{5} \sqrt[3]{2} 
              - \frac{1}{5} \sqrt[3]{4}\quad \text{and} \quad
   \alpha_2 = -\frac{16641}{7660} + \frac{1290}{383} \sqrt[3]{2} 
                   + \frac{1000}{383} \sqrt[3]{4} .  $$
We find
\begin{align*}
  s_1 & =  - \frac35             & t_1 & =  - \frac15    & 
      S_- & = - \frac{342}{383}  & T_+ & = - \frac{200}{383}  \\
  s_2 & = \frac{1290}{383}       & t_2 & = \frac{1000}{383} & 
      S_+ & = - \frac{858}{383}  & T_- & = - \frac{5383}{1915}  \\
  \Sigma & =  -\frac{6138414}{733445}  &  &  &  &    
\end{align*}
This yields 
$$ s_3 = - \frac{28099233}{66234835}, \qquad
   t_3 = - \frac{5000211}{66234835}, \qquad 
   r_3 = - \frac{s_3^2}{2t_3}
       = \frac{27002048329}{22652313570} $$
in perfect agreement with our calculations at the end of Section \ref{S1}.

\section{Squares of the Form $a - b \omega$}

Consider more generally the problem of classifying squares of the form
$a - b \omega$ for some fixed value of $b \in \Q^\times$. Since
$a - b \omega = a - \omega'$ for $\omega' = \sqrt[3]{mb^3}$ we find,
by simply replacing $m$ with $mb^3$ in Thm. \ref{Tm}, the following

\begin{thm}
Let $m$ be an integer that is not a cube, let $K = \Q(\omega)$ with 
$\omega^3 = m$ denote the corresponding pure cubic field, and let 
$E_m$ denote the elliptic curve $y^2 = x^3 - m$. There is a bijection 
between elements $\alpha  \in K^\times$ with $\alpha^2 = a - b\omega$ 
for $a, b \in \Q$ and the rational points $(x,y) \in E_m'(\Q)$,
where $E_m': y^2 = x^3 - mb^3$ is a quadratic twist of $E_m$. 

In fact, if $\alpha = r + s \omega + t \omega^2$ satisfies 
$\alpha^2 = a - b\omega$, then $(x,y) \in E_m'(\Q)$ for 
$x = bs/t$ and $y = b^2/t$. 

Conversely, an affine point $(x,y) \in E_m'(\Q)$ corresponds to 
$$ \alpha = -\frac{x^2}{2y} + \frac {bx}{y} \omega + \frac{b^2}y \omega^2 
     \quad \text{ with } \alpha^2 = \frac{x^4 + 8mx}{4y^2}  - b \omega. $$

In particular, $a - b\omega$ is a square with $N(a - b\omega) = y^2$ if 
and only if $(a,y) \in 2 E'(\Q)$. 
\end{thm}

\medskip \noindent{\bf Example.} 
The unit
$$ -19 + 7 \sqrt[3]{20} = (1 + \sqrt[3]{20} - \sqrt[3]{50}\,)^2 $$
is a well known example due to Nagell \cite{Nag}; since
$\sqrt[3]{50} = \frac12\sqrt[3]{20}^2$, it  comes from the point 
$P(-2,-2)$ on the quadratic twist $-7y^2 = x^3 - 20$ of the
elliptic curve $y^2 = x^3 - 20$. The corresponding point
on the elliptic curve $y^2 = x^3 + 7^3 \cdot 20$ is $P'(14,98)$.
Observe that $2P' = (-19,1)$.

\begin{cor}
There is a binomial unit $a - b \omega$ that is a square in $K$
if and only of there is an integral point $(a,1) \in 2E(\Q)$ on
$E: y^2 = x^3 - mb^3$.
\end{cor}

\section{A Homomorphism from $E(\Q)$ to $\Cl(K)$}

In this section we will define a map from $E(\Q)$, the group
of rational points on $E: y^2 = x^3 - m$ (where
$m \not \equiv 0, \pm 1 \bmod 9$ is a cubefree integer), to the 
group $\Cl(K)[2]$ of ideal classes of order dividing $2$ in 
$K = \Q(\sqrt[3]{m}\,)$, and show that this map is a homomorphism.

\begin{thm}\label{Tids}
Let $m \not\equiv 0, \pm 1 \bmod 9$ be a cubefree integer, and
consider the cubic number field $K = \Q(\sqrt[3]{m}\,)$ and
the elliptic curve $E: y^2 = x^3 - m$. For any 
$P = (x,y) \in E(\Q) \setminus \{\infty\}$, set
$x = a/e^2$ for coprime integers $a$ and $e$. Then
$(a - e^2 \omega) = \fa^2$, and the map $P \to [\fa]$ induces a 
homomorphism $\kappa: E(\Q) \lra \Cl(K)[2]$.
\end{thm}

\begin{proof}
Since $a^3/b^3 - m = y^2$ we find that $N(a - b \omega) = y^2b^3$.
Since $b = e^2$ is a square (see \cite[p. 68]{ST}), the norm of 
$\alpha = a - b \omega$ is a square. Thus $(\alpha)$ is an
ideal square if and only if $a-b\omega$ is coprime to 
its conjugates $\alpha' = a - b \omega \rho$ and 
$\alpha'' = a - b \omega \rho^2$ in the normal closure $L$ of $K/\Q$.
Let $\fd$ be the greatest common ideal divisor of $\alpha$
and $\alpha'$. Then $\fd  \mid \alpha''$, hence
$\fd$ divides the trace $3a$ of $\alpha$ as well as the difference
$\alpha - \alpha' = b\omega(1-\rho)$. Since $a$ and $b$ are coprime,
$\fd$ is a product of ideals above $3$.  

Since we have assumed that $m \not\equiv \pm 1 \bmod 8$, we have
$(3) = \fq^3$ in $K$ (see e.g. \cite{Cohn} for the decomposition
law in pure cubic number fields). Assume therefore that $\fq \mid (\alpha)$
(this implies that $\fq \mid (\alpha')$ since $\fq$ is totally 
ramified). If $3 \mid m$, then $3 \mid a$, and from 
$y^2 = a^3/b^3 - m$ we deduce that $3^2 \mid m$ contradicting our
assumptions. If $3 \nmid m$, then $a^3 - mb^3 = y^2 b^3$ is divisible
by $3$ if and only if it is divisibly by $9$, and now
$a^3 \equiv mb^3 \bmod 9$ implies that $m \equiv \pm 1 \bmod 9$,
again contradicting our assumptions. Thus $(\alpha) = \fa^2$, and 
$[\fa] = \kappa(P)$ is an ideal class of order dividing $2$ in $K$.
It remains to show that $\kappa$ is a homomorphism.

To this end, assume that $x_P = a/b$ and $x_Q = c/d$ for $P, Q \in E(\Q)$. 
Set $\kappa(P) = [\fa]$ and $\kappa(Q) = [\fb]$. We have to show
that $\kappa(P+Q) = [\fa\fb]$, which is equivalent to 
$\kappa(P) \kappa(Q) \kappa(R) \sim 1$ for collinear rational 
points  $P, Q, R \in E(\Q)$.

The Weil map $\alpha: E(K) \lra K^\times/\Kts$ defined by
$$ \alpha(P) = \begin{cases}
            (x - \omega) \Kts & \text{ if } P = (x,y), \ y \ne 0 \\
             3 \omega    \Kts & \text{ if } P = (\omega,0),  \\ 
                         \Kts & \text{ if } P = \cO  
               \end{cases} $$
(see e.g. \cite[(4.2)]{Huse}) is a homomorphism, and if we restrict
$\alpha$ to $E(\Q)$, then the second case cannot occur. This shows that
$(x_P - \omega)(x_Q - \omega)(x_R - \omega)$ is a square in $K^\times$.
Thus $(a - b \omega)(c - d \omega)(e - f\omega) = \beta^2$, where
$x_R = e/f$, and $\fa^2 \fb^2 \fc^2 = (\beta)^2$. But this implies 
$\fa \fb \fc = (\beta)$.
\end{proof}

We now discuss a few examples.

\medskip\noindent{\bf Example 1.}
Let $K = \Q(\sqrt[3]{26}\,)$ and $E: y^2 = x^3 - 26$. 
Then $E(\Q)$ is generated by $P=(3,1)$ and $Q = (35,207)$.
The element $3 - \omega$ is a unit, and 
$(35 - \omega) = \fthr_1^3 \fthr_2 \fq^2$, where
$(3) = \fthr_1 \fthr_2^2$ and where $\fq$ is a prime ideal
above $23$.

This shows that the condition $m \not \equiv \pm 1 \bmod 9$
cannot be dropped.

\medskip\noindent{\bf Example 2.}
Let $K = \Q(\sqrt[3]{47}\,)$ and $E: y^2 = x^3 - 47$. 
Then $E(\Q)$ is generated by $P=(6,13)$ and $Q = (12,41)$;
moreover, $P+Q = (34/9, -71/27)$. We have
$$  (6 -  \omega) = \fp_{13}^2, \quad
 (12 -  \omega)   = \fp_{41}^2, \quad
 (34 - 9\omega)   = \fp_{71}^2. $$

The ideals $\fp_{13}$ and $\fp_{41}$ generate the ideal class
of order $2$, the ideal $\fp_{71}$ is principal.

\medskip\noindent
{\bf Example 3.}
Let $m = 57$; then $E(\Q)$ is generated by $P = (4873/36, -340165/216)$,
$(4873 - 36 \omega) = \fa^2$, and $\fa$ is principal. The elliptic 
curve $E': y^2 = x^3 - 57^2$ has rank $0$ according to {\tt sage} since 
$L(E',1) \ne 0$, and has nontrivial Tate-Shafarevich group $\TS(E')$. 

This suggests the question whether there is any connection between
the kernel of $\kappa$ in $E(\Q)$ and the Tate-Shafarevich group 
$\TS(E')$ of $E'$.

\medskip\noindent
{\bf Remark.} If we consider the more general case of elliptic
curves $y^2 = x^3 - mb^3$, then a point $(x,y) \in E(\Q)$ with 
$x = a/e^2$ provides us with an element $ \alpha = a-be^2\omega$
with square norm $N(\alpha) = a^3 - mb^3e^6 = (ye^3)^2$; yet
$(\alpha)$ need not be an ideal square if $\gcd(a,b) \ne 1$.
In fact, $(14,44)$ is a rational point on $y^2 = x^3 - 2^3 \cdot 101$,
and neither $(14 - 2 \omega)$ nor $(7 - \omega)$ are ideal squares
in $K = \Q(\omega)$ with $\omega^3 = 101$.

On the other hand, the proof of Thm. \ref{Tids} shows

\begin{lem}\label{Lids}
Assume that $m$ is a cubefree integer with $m \not\equiv 0, \pm 1 \bmod 9$,
set $\omega = \sqrt[3]{m}$, and let $P = (x,y)$ be a rational point on 
$E: y^2 = x^3 -mb^3$ for some integer $b$. If $x = a/e^2$, then 
$(a - be^2\omega)$ is an ideal square in $K$ if $\gcd(a,b) = 1$.
\end{lem}

\section{Hilbert Class Fields}

Let $K = \Q(\omega)$ with $\omega^3 = m$ be a cubic number 
field, where we assume that $m \not\equiv 0, \pm 1 \bmod 9$, and 
let $E: y^2 = x^3 - mb^3$ be a quadratic twist of $E$.
If $P = (x,y) \in E(\Q)\setminus 2E(\Q)$ is a rational point
with $x = r/s^2$ and coprime integers $r, s$, then 
$\alpha = r - s^2b\omega \in K$ has norm 
$r^3 - ms^6b^6 = y^2(bs)^6$. Since $\alpha$ is not a square,
$K_P = K(\sqrt{\alpha}\,)$ is a quadratic extension; since 
$(\alpha)$ is an ideal square by Lemma \ref{Lids}, $L/K$ can only
ramify at the primes above $2$ and $\infty$. If $2 \mid s$ and if 
$r$ is chosen positive, then $r \equiv 1 \bmod 4$ since 
$r^3 \equiv 1 \bmod 4$, and in this case the extension $K_P$ is 
unramified everywhere.
 
Table \ref{T1} lists the pure cubic number fields with
even class number and $m \le 113$, an elliptic curve $E$ whose
rational points provide us with unramified quadratic extensions
$K(\sqrt{\alpha}\,)$.

\begin{table}[h!]
$$ \begin{array}{r|clcc}
      m  &  \Cl_2(K) &   E   &     x_P    &  \alpha \\  \hline
     11  & 2  & y^2 = x^3 -  11 &   9/4    &  9 - 4 \omega \\ 
     15  & 2  & y^2 = x^3 +  15 &   1/4    &  1 + 4 \omega \\
     39  & 2  & y^2 = x^3 +  39 & 217/4    & 217 + 4\omega \\
     43  & 4  & y^2 = x^3 -  43 & 1177/36  & 1177 - 36\omega \\
     47  & 2  & y^2 = x^3 +  47 &   17/4   &  17 - 4 \omega \\
     57  & 2  & y^2 = x^3 -  57 & 4873/36  & 4873 - 36 \omega \\
     58  & 2  & y^2 = x^3 -  58 & 5393/484 & 5393 - 484 \omega \\
     61  & 2  & y^2 = x^3 -  61 & 929/100  & 929 - 100 \omega \\
     63  & 2  & y^2 = x^3 +  63 &    9/4   &  9 + 4 \omega \\
     65  & 2  & y^2 = x^3 +  27 \cdot 65   & 129/4 & 129 + 3 \cdot 4 \omega \\
     66  & 2  & y^2 = x^3 +  66 &    1/4   &  1 + 4 \omega \\
     67  & 2  & y^2 = x^3 -  67 &   17/4   & 17 - 4 \omega \\
     76  & 2  & y^2 = x^3 -  76 &   17/4   & 17 - 4 \omega \\
     79  & 2  & y^2 = x^3 +  27 \cdot 79   
                            & 1921/100 & 1921 + 3 \cdot 100 \omega \\
     83  & 2  & y^2 = x^3 -  83 &   33/4   &  33 - 4 \omega \\
     89  & 2  & y^2 = x^3 -  89 &  153/4   & 153 - 4 \omega \\ 
    101  & 2  & y^2 = x^3 - 101 & 6342921/1073296 & 6342921 - 1073296 \omega \\
    105  & 2  & y^2 = x^3 - 105^2 & 16465/196 & 16465 - 196 \omega^2 \\
    106  & 2  & y^2 = x^3 - 106   & 8297/1024 & 8297 - 1024 \omega \\
    113  & (2,2) & y^2 = x^3 - 27 \cdot 113 &  97/4 & 97 - 3 \cdot 4 \omega \\
         &       &             & 43449/2500  & 43449 - 3 \cdot 2500 \omega\\
         &       &             & 1257/64  & 1257- 3 \cdot 64 \omega  \\
   2351  & (4,2,2) & y^2 = x^3 + 27m 
                           &   57/4  &   57 + 3 \cdot 4 \omega \\
         &         &       & -551/16 & -551 + 3 \cdot 16 \omega \\
         &         &       & - 87/4  & - 87 + 3 \cdot 4 \omega                
   \end{array} $$
\caption{Quadratic Unramified Extensions of 
                   Pure Cubic Number Fields}\label{T1}
\end{table}

If $m = 113$, the three quadratic unramified extensions are
\begin{align*}
  f_1(x) & = x^6 - 291x^4 + 28227x^2 - 717409 \\
  f_2(x) & = x^6 - 130347x^4 + 5663446803x^2 - 34351825047849 \\ 
  f_3(x) & = x^6 - 3771x^4 + 4740147x^2 - 1186320249
\end{align*}
Since neither of the three points from which these extensions
originate are in $2E(\Q)$, these three extensions are pairwise
distinct.

In the case of $m = 2351$, it can be checked that the three 
unramified (in all three cases, we have $\alpha  > 0$, so there
is no ramification at infinity) quadratic extensions generated
by a root of the polynomials
 \begin{align*}
  f_1(x) & = x^6 -  171x^4 +   9747x^2 - 4247721, \\
  f_2(x) & = x^6 + 1653x^4 + 910803x^2 - 92717641, \\
  f_3(x) & = x^6 +  261x^4 +  22707x^2 - 3404025
\end{align*}
are independent.

Our results suggest that every unramified quadratic extension
of a pure cubic number field $K = \Q(\sqrt[3]{m}\,)$ can be 
computed from a rational point on some quadratic twist of the
elliptic curve $y^2 = x^3 - m$.

\vskip 0.5cm

\end{document}